\documentclass[12pt]{article} \textheight 8.8 true in \textwidth
6.2 true in

\usepackage[T2A]{fontenc}
\usepackage[cp1251]{inputenc}
\usepackage[english]{babel}
\usepackage{amsfonts}
\usepackage{amssymb,amsmath,amstext}
\usepackage{amscd, amsthm, latexsym}
\usepackage{color}
\usepackage[final]{graphicx}
\usepackage{epsfig}
\usepackage{euscript}
\usepackage{srcltx,epsf}
\usepackage{wrapfig}
\usepackage{mathrsfs}
\usepackage{ifthen}
\usepackage{cleveref}

\def\hide#1{}
\def\old#1{}

\def\oop#1{}
\def\gap#1{}

\newcommand{\F}[1]{\mathfrak #1}

\theoremstyle{plain}
\newtheorem{theorem}{Theorem}
\newtheorem*{theorem*}{Theorem}
\newtheorem{lemma}{Lemma}
\newtheorem{remark}{Remark}

\newtheorem{definition}{Definition}

\def \ZZ  {\Bbb Z}

\def\RR{\mathbb R}
\begin{document}

\newcounter{figcounter}
\setcounter{figcounter}{0} \addtocounter{figcounter}{1}

\title{On the invariants of $4\times 4$ skew-symmetric matrix under cluster mutations.}
\author{G.~Chelnokov {\small\tt grishabenruven@yandex.ru }
\\ {\small\em National Research University Higher School of Economics, Moscow, Russia} \\[2ex]
 }
\date{}
\maketitle \maketitle

\begin{abstract}
We prove that any invariant of a 4-quiver, that is piecewise
polynomial, moreover, polynomial for fixed signs of entries, is a
function of determinant of a quiver.


\end{abstract}

The main purpose of this paper is to prove Theorem \ref{only_Det},
but to formulate it smoothly we need some definitions.

\begin{definition}
Call a {\em quiver of size} $n$ (or {\em $n$-quiver}) an $n\times n$
skew-symmetric matrix with real entries. By $U_n$ denote the set of
all $n$-quivers. Call {\em a carriage} $S_i$ a subset of $U_n$ of
all quivers with prescribed signs of entries, so $U_n$ is covered by
$2^{\binom{n}{2}}$ carriages. Call a quiver $X$ {\em an inner point}
of its carriage if all entries of $X$ are nonzero.
\end{definition}

\begin{definition}
Given an integer $k \in [1..n]$ by $\mu_k$ denote the mapping $U_n
\to U_n$, such that under the mapping $X \to \mu_k(X)$ $x_{i,j} \to
-x_{i,j}$ if $k=i\; \text{or}\; j$. If $k\neq i,j$ then $x_{i,j} \to
x_{i,j}+x_{i,k}x_{k,j}$ if both $x_{i,k}$ and $x_{k,j}$ are
positive,  $x_{i,j} \to x_{i,j}-x_{i,k}x_{k,j}$ if both $x_{i,k}$
and $x_{k,j}$ are negative, else  $x_{i,j} \to x_{i,j}$. Call
$\mu_k$ {\em a cluster mutation with respect to vertex $k$}.
\end{definition}

\begin{definition}
Let $\mathcal{P}$ be an $2^{\binom{n}{2}}$-plet of polynomials $P_i$
in matrix entries (we will write $P_i(X)$ instead of
$P_i(x_{1,2},\dots,x_{n-1,n})$, keep in mind that it is not a
polynomial in matrix sense), moreover, let polynomials $P_i$
bijectively correspond to carriages $S_i$. By $F_\mathcal{P}: U_n
\to \RR$ denote the function $F_\mathcal{P}: X\in S_i \to P_i(X)$.
We call such function {\em carriages-wise polynomial}.
\end{definition}

Certainly we say that $F_\mathcal{P}$ is invariant under cluster
mutations if $F_\mathcal{P}(X)=F_\mathcal{P}(\mu_k(X))$ for all $k$
and all $X\in U_n$. Now we have all the notions necessary to
formulate the main result of this paper.

\begin{theorem}\label{only_Det}
Let $F_\mathcal{P}$ be an carriages-wise polynomial function,
invariant under cluster mutations, on 4-quivers.
Then there exists a polynomial in one variable $f$ such that
$F_\mathcal{P}(X)=f(Det(X))$.
\end{theorem}

\begin{proof} Essentially in this prove we will repeatedly use
an obvious fact that if two (multivariable) polynomials coincide as
functions on an open set then they coincide as polynomials. First we
need one more auxiliary notion.

\begin{definition}
Let $S_i$ be a carriage in $U_n$. By $\mu\mu_k^i$ denote a mapping
$\mu\mu_k^i: U_n \to U_n$, which is polynomial in all entries and
coincide with $\mu_k(X)$ for $X\in S_i$. For example if $\mu_k(X)$
maps $x_{1,2} \to x_{1,2} + x_{1,k}x_{k,2}$ for $X\in S_i$ then
 $\mu\mu_k^i$ maps $x_{1,2} \to x_{1,2} +
x_{1,k}x_{k,2}$ for all $X$.
\end{definition}

\begin{remark}{\label{rem1}} If two carriages $S_{i_1}$ and $S_{i_2}$ differ only in the orientation of some arrows not adjusted to
vertex $k$, then the mappings $\mu\mu_k^{i_1}$ and $\mu\mu_k^{i_2}$
coincide.
\end{remark}

\begin{lemma}\label{through_one_edge}
Let $F_\mathcal{P}$ be an carriages-wise polynomial function,
invariant under cluster mutations, on $n$-quivers. Assume two
carriages $S_{i_1}$ and $S_{i_2}$ differ in a sign of just one entry
$x_{i,j}$, and there exist a vertex $k$, such that $x_{i,k}$ and
$x_{k,j}$ have same sings. Then $P_{i_1}$ and $P_{i_2}$ coincide.
\end{lemma}

\begin{proof}
Consider inner quivers $X_1\in S_{i_1}$ and $X_2\in S_{i_2}$, such
that all their entries other then $x_{i,j}$ coincide, entry
$x_{i,j}$ differ only in sign and $|x_{i,j}|<x_{i,k}x_{k,j}$. Then
$\mu_k$ maps $X_1$ and $X_2$ into inner quivers of the same
carriage, denote this carriage $S_j$.

Polynomial mappings $\mu\mu_k^{i_1}$ and $\mu\mu_k^{i_2}$ coincide
due to Remark \ref{rem1}. Also $\mu\mu_k^{j} \circ \mu\mu_k^{i_1}
=id$ because $\mu_k(\mu\mu_k^{i_1}(X))=X$ holds for all $X\in
S_{i_1} \cap \mu_k(S_j)$ (indeed, here we use the fact that $\mu_k$
is an involution, otherwise we would have to write $X \in S_{i_1}
\cap \mu_k^{-1}(S_j)$), and $\mu\mu_k^{j} \circ \mu\mu_k^{i_1}$ is a
polynomial mapping.

Since $F_\mathcal{P}$ is invariant, $P_{i_1}(X)=P_j(\mu_k(X))$ for
all $X \in S_{i_1} \cap \mu_k(S_j)$. Then
$P_{i_1}(X)=P_j(\mu\mu_k^{i_1}(X))$ holds as polynomial equality,
that is for arbitrary $X$. Similarly we obtain
$P_{j}(X_1)=P_{i_2}(\mu\mu_k^j(X_1))$.

Substituting into the latter one $X_1=\mu\mu_k^{i_1}(X))$, combining
with the first one and applying $\mu\mu_k^{j} \circ \mu\mu_k^{i_1}
=id$ we get $P_{i_1}(X)=P_{i_2}(X)$ as desired.
\end{proof}

\begin{lemma}\label{easy_combi}
The statement of the Lemma \ref{through_one_edge} implies that in
case $n=4$ all $P_i$ coincide.
\end{lemma}

\begin{proof} Let us reformulate our statement. Consider four
vertices, each pair is connected by an arrow in one of two
directions. We are  allowed for an oriented path
$\overrightarrow{AB}, \overrightarrow{BC}$ to switch the direction
of an arrow between $A$ and $C$. We need to prove that by means of
such operations we can transit from any configuration to any one.
Here a configuration of arrows represent a carriage, and our right
to pass from configuration $S_i$ to $S_j$ represent that Lemma
\ref{through_one_edge} claims $P_i=P_j$.

The proof of the above statement consists of a  cases consideration.
Call a vertex {\em regular} in some configuration if this vertex
have positive indegree and outdegree. First, if $A$ is regular then
we can switch arrows between $B,C,D$ as we want. Indeed, without
loss of generality assume arrows are $\overrightarrow{BA}$,
$\overrightarrow{AC}$ and $\overrightarrow{AD}$; then we can switch
$(B,C)$ and $(B,D)$ as we want; then we can switch $(B,C)$ and
$(B,D)$ to $\overrightarrow{BC}$ and $\overrightarrow{CD}$, then
switch $(B,D)$ as we want, then switch $(B,C)$ and $(B,D)$ as we
want. Second, we can pass between any to configurations such that
$A$ is regular. Third, each configuration have a regular vertex.
Assume initial configuration have a regular vertex $A$ and final one
have a regular vertex $B$. Then using first statement we can make
vertex $B$ regular, then using second statement bring adjusted to
$B$ arrows into the correct position, then again use the first
statement.
\end{proof}

\begin{remark} Lemma \ref{easy_combi} is an essential step that do
not works in case of $3\times 3$ matrices. Denote the entries by
$$
\begin{pmatrix}  0 & x & -y \\
-x & 0 & z \\ y & -z & 0\end{pmatrix} $$ Then Lemma
\ref{through_one_edge} claims $P_{x\geqslant 0, y\geqslant 0,
z\geqslant 0}=P_{x\leqslant 0, y\geqslant 0, z\geqslant
0}=P_{x\geqslant 0, y\leqslant 0, z\geqslant 0}=P_{x\geqslant 0,
y\geqslant 0, z\leqslant 0}$ and $P_{x\leqslant 0, y\leqslant 0,
z\leqslant 0}=P_{x\geqslant 0, y\leqslant 0, z\leqslant
0}=P_{x\leqslant 0, y\geqslant 0, z\leqslant 0}=P_{x\leqslant 0,
y\leqslant 0, z\geqslant 0}$, but nothing else. Thanks to this
becomes possible the invariant
$$
F_{\mathcal{P}}(X)=\left\{\begin{aligned}
x^2+y^2+z^2-xyz \quad \text{if 3 or 2 non-negative among x,y,z}  \\
x^2+y^2+z^2+xyz \quad \text{if 3 or 2 non-positive among x,y,z}  \\
\end{aligned} \right.
$$
\end{remark}

To formulate the next lemma smoothly we introduce some new
notations. We address entries of a skew-symmetric matrix as
$$
X=\begin{pmatrix}  0 & x & y & z \\
-x & 0 & u & -v \\ -y & -u & 0 & w \\ -z & v & -w & 0\end{pmatrix}
$$ First, by introducing different letters we avoid index
pandemonium, second, this is a more symmetric orientation (if
$x,y,z,u,v,w >0$ arrows are oriented $\overrightarrow{32}$,
$\overrightarrow{24}$, $\overrightarrow{43}$ and
$\overrightarrow{i1}$). In this notation $Det(X)=(xw+yv+zu)^2$.

Next, we address vector $(x,y,z)$ as $Y$ and vector $(u,v,w)$ as
$V$. So, short for $F_{\mathcal{P}}(x,y,z,u,v,w)$ is
$F_{\mathcal{P}}(Y,V)$. By $(\cdot)$ we denote bilinear form $Y\cdot
V=xw+yv+zu$.

\begin{lemma}\label{to scolar prod}
Let $F_\mathcal{P}$ be an carriages-wise polynomial function,
invariant under cluster mutations, on 4-quivers. Then
$F_\mathcal{P}(Y,V)=F_\mathcal{P}(Y+Y_1,V)$ holds for any $Y_1$ such
that $Y_1\cdot V=0$, and $F_\mathcal{P}(Y,V)=F_\mathcal{P}(Y,V+V_1)$
holds for any $V_1$ such that $Y\cdot V_1=0$.
\end{lemma}

\begin{proof}
By virtue of Lemmas \ref{through_one_edge},\ref{easy_combi} we
consider $F_\mathcal{P}$ to be polynomial.

Consider a carriage given by $v,z \geqslant 0$ and $w\leqslant 0$,
other three signs unimportant. Denote this carriage $S_1$. Then the
explicit formula for $\mu\mu_4^1$ is
$\mu\mu_4^1(x,y,z,u,v,w)=(x+vz,y+(-w)z,-z,u,-v,-w)$. So
$\big(\mu\mu_4^1\big)^2(x,y,z,u,v,w)=(x+2vz,y-2wz,z,u,v,w)$ or in
vector notation
$\big(\mu\mu_4^1\big)^2(x,y,z,u,v,w)=(x,y,z,u,v,w)+2z(v,-w,0,0,0,0)$.
Consequently
$\big(\mu\mu_4^1\big)^{2k}(x,y,z,u,v,w)=(x,y,z,u,v,w)+2kz(v,-w,0,0,0,0)$
for all integer $k$. So the restriction of polynomial
$F_\mathcal{P}$ to the line
$\{(x,y,z,u,v,w)+\lambda(v,-w,0,0,0,0)|\; \lambda \in \RR\}$ takes
same value infinitely many times, thus this restriction is a
constant. Similarly, by considering carriage $w,y \geqslant
0,\;\;u\leqslant 0$ and vertex 3 get that $F_\mathcal{P}$ is
constant on the line $\{(x,y,z,u,v,w)+\lambda(-u,0,w,0,0,0)|\;
\lambda \in \RR\}$. Since any vector $Y_1$ such that $Y_1\cdot V=0$
belongs to the linear hull $\langle (v,-w,0,0,0,0),
(-u,0,w,0,0,0)\rangle$ we are done with the statement
$F_\mathcal{P}(Y,V)=F_\mathcal{P}(Y+Y_1,V)$.

The second statement of the Lemma can be achieved similarly.
\end{proof}

\begin{lemma}\label{final}
Let $F$ be an function $F: \RR^3\times\RR^3 \to \RR$ satisfying
conditions \begin{itemize}
\item $F(Y_1,V)=F(Y_2,V)$ holds for any
$Y_1,Y_2,V$ such that $Y_1\cdot V=Y_2\cdot V$,
\item $F(Y,V_1)=F(Y,V_2)$ holds for
any $V_1,V_2,Y$ such that $Y\cdot V_1=Y\cdot V_2$. \end{itemize}
Then there exist a function $f: \RR \to \RR$, such that
$F(Y,V)=f(Y\cdot V)$.
\end{lemma}

Note that the premise of Lemma \ref{final} is exactly the conclusion
of Lemma \ref{to scolar prod} formulated in more symmetric terms.

\begin{proof}
Define $f(x)=F((x,0,0),(1,0,0))$. Our goal is to transit from
$F((y_1,y_2,y_3),(v_1,v_2,v_3))$ via a sequence of equalities,
provided by Lemma \ref{final}. First, we may assume one of $y_1,v_1$
be nonzero, otherwise
$F((0,y_2,y_3),(0,v_2,v_3))=F((0,y_2,y_3),(1,v_2,v_3))$. Without
loss of generality $v_1\neq 0$. Then

\begin{align*}
F((y_1,y_2,y_3),(v_1,v_2,v_3))=F(\frac{y_1v_1+y_2v_2+y_3v_3}{v_1},0,0),(v_1,v_2,v_3)=\\F(\frac{y_1v_1+y_2v_2+y_3v_3}{v_1},0,0),(v_1,1,0)=
F((0,y_1v_1+y_2v_2+y_3v_3,0),(v_1,1,0))=\\F((0,y_1v_1+y_2v_2+y_3v_3,0),(1,1,0))=F((y_1v_1+y_2v_2+y_3v_3,0,0),(1,1,0))=\\F((y_1v_1+y_2v_2+y_3v_3,0,0),(1,0,0)).
\end{align*}
\end{proof}
Lemma \ref{final} finishes the proof of Theorem \ref{only_Det}.
\end{proof}

\subsection*{Some final remarks}

Indeed, analogues of Theorem \ref{only_Det} can be proved  for wider
classes of functions by the arguments of this paper. We use just two
properties of the class of polynomials. First, that equality of two
functions on an open set implies global equality. This holds for
analytic functions. Second, that a function, taking on a
1-dimensional affine subspace a fixed value infinitely many times is
a constant on this line. This holds (for example) for rational
functions and real exponents, but not complex exponents. So, the
result follows for rational functions, real quasi-polynomials, real
quasi-rational similarly.

Also, the question of invariants of integer quivers is more widely
known then it's reals counterpart. In this case the difference is
insignificant. Indeed, the intersections of (real)carriages with
integer lattice $\ZZ^6$ are sufficiently ``large'' to any
carriage-wise polynomial (rational, real quasi-rational...)
invariant on integer quivers be an invariant on real quivers.

\end{document}